\documentclass{amsart}
\usepackage{amsmath,amsthm,amsfonts,amssymb}
\usepackage{graphicx}
\usepackage{psfrag}
\newtheorem{lem}{Lemma}
\newtheorem{prop}{Proposition}
\newtheorem{thm}{Theorem}
\newtheorem*{defn}{Definition}
\newtheorem{cor}{Corollary}

\newcommand{\N}{\mathbb{N}} \newcommand{\T}{\mathbb{T}}
\newcommand{\Q}{\mathbb{Q}} \newcommand{\R}{\mathbb{R}}
\newcommand{\Z}{\mathbb{Z}} 

\newcommand{\tnorm}[1]{|\!|\!| #1 |\!|\!|}

\DeclareMathOperator{\diff}{Diff}
\DeclareMathOperator{\diam}{diam}
\begin{document}

\title{Nonstandard Smooth Realizations of Liouville Rotations}

\author[B. R. Fayad]{Bassam Fayad} 

\address{Bassam Fayad, LAGA, UMR 7539,
  Universit\'e Paris 13, 93430 Villetaneuse, France}

\author[M. Saprykina]{Maria Saprykina} 

\address{Maria Saprykina, Department of Math \& Stats, Jeffery Hall,
  University Ave.  Kingston, ON Canada, K7L 3N6 }

\thanks{M. Saprykina would like to thank the University of Texas at
  Austin and Rafael de la Llave for support and hospitality.}
  
\author[A. Windsor]{Alistair Windsor} 

\address{Alistair Windsor, Department of Mathematics, University of
  Texas at Austin, 1 University Station, C1200, Austin, TX 78712-0257,
  USA }

\thanks{A. Windsor gratefully acknowledges the generous support of the
  Universit\'e Paris 13 and the CNRS as well as the hospitality of the
  Laboratoire de Probabilit\'es et Mod\`eles Al\'eatoires.}

\begin{abstract}
  We augment the method of $C^\infty$ conjugation approximation with
  explicit estimates on the conjugacy map. This allows us to construct
  ergodic volume preserving diffeomorphisms measure-theoretically
  isomorphic to any apriori given Liouville rotation on a variety of
  manifolds. In the special case of tori the maps can be made uniquely
  ergodic.
\end{abstract}

\maketitle 
\section{Introduction}
\label{sec:Introduction}
We call a diffeomorphism $f$ of a compact manifold $M$ that preserves
a smooth measure $\mu$ a \emph{smooth realization} of an abstract
system $(X,T,\nu)$ if they are measure-theoretically isomorphic. A
diffeomorphism of a compact manifold has finite entropy with respect
to any Borel measure. The natural question therefore becomes whether
every finite entropy automorphism of a Lebesgue space has a smooth
realization. This problem remain stubbornly intractable and there
remain abstract examples that have no known smooth realizations.

We seek to find smooth realizations of one of the simplest types of
automorphisms; aperiodic automorphisms with pure point spectrum with a
group of eignevalues with a single generator. Such automorphisms are
measure theoretically isomorphic to irrational rotations of the
circle. They therefore have a natural smooth realization. We seek
smooth realizations on manifolds other than $\T$. Such relizations are
called non-standard smooth realizations. 

We extend the conjugation approximation method of Anosov and Katok
\cite{AK} to construct non-standard smooth realizations of a given
Liouville rotation on $\T$ on a variety of manifolds $M$. Indeed, in
the special case that the manifold is $\T^d$ for $d \geq 2$, we can
produce uniquely ergodic realizations of the given Liouville
rotation. The crucial new ingredient is an explicit construction of
the conjugating maps that allows us to estimate their
derivatives. This allows us to ensure that the construction converges
for a predetermined Liouville number $\alpha$. The approach parallels
that taken in \cite{FS}.  The original paper of Anosov and Katok paper
constructed non-standard smooth realizations of a dense set of
Liouville rotations. However, without estimates, it was not possible
to identify which Liouville rotations could be realized.

\begin{defn}
  A number $\alpha \in \R \backslash \Q$ is a Liouville number if for
  all $k >0$ we have
  \begin{equation}
    \label{eq:Liouville}
    \liminf_{q \rightarrow \infty} q^k \| q \alpha \| = 0
  \end{equation}
  where $\|q \alpha \| = \inf_{p \in \Z} | q \alpha - p|$.
\end{defn}

Let $\T^d:= \R^d / \Z^d$ denote the $d$-dimensional torus. Let
$R_\theta : \T \rightarrow \T$ be the rotation of the circle, taken
with the Haar probability measure, given by
$R_\theta(x) = x + \theta \mod 1$. 

Denote by $\diff^{\infty}( M , \mu)$ the class of $C^\infty$
diffeomorphisms of $M$ that preserve a $C^\infty$  smooth volume
$\mu$. Throughout this paper we will use $\lambda$ for the probability
measure induced by the standard Lebesgue measure.

\begin{thm}\label{thm:GeneralCase}
  Let $M$ be a compact connected manifold, possibly with boundary, of
  dimension at least 2 that admits an effective $C^\infty$ action of
  $\T$ preserving a $C^\infty$ smooth volume $\mu$. For every $\alpha
  \in \R \backslash \Q$ Liouville there exists an ergodic $T \in
  \diff^{\infty}(M, \mu)$ measure-theoretically isomorphic to the
  rotation $R_\alpha$.
\end{thm}

In the special case $M = \T^d$ we can strengthen the result to obtain
unique ergodicity. 

\begin{thm}\label{thm:TorusCase}
  For every Liouville $\alpha \in \R \backslash \Q$, and every $d \geq
  2$ there exists a uniquely ergodic transformation $T \in
  \diff^{\infty}(\T^d, \lambda)$ such that $T$ is
  measure-theoretically isomorphic to the rotation $R_\alpha$.
\end{thm}

It remains open whether there are $C^\infty$ realizations of Diophantine
rotations on any manifold other than $\T$. 

\section{Construction}
\label{sec:Construction}
\subsection{Outline}

The required measure preserving diffeomorphism $T$ is constructed as
the limit of a sequence of periodic measure preserving diffeomorphisms
$T_n$. For each of the properties that we wish the limiting
diffeomorphism $T$ to possess, we establish an appropriate finitary
version possessed by the periodic diffeomorphism $T_n$.

Let $S:\T \times M \rightarrow M$ denote an effective $C^\infty$
action of $\T$ on $M$ that preserves the volume and denote by
$S_\alpha$ the diffeomorphism $S(\alpha, \cdot)$. The diffeomorphism
$T_n$ is given by
\begin{equation}
  \label{eq:TnDefn}
  T_n := H_n S_{\alpha_n} H_n^{-1}
\end{equation}
where $\alpha_n \in \Q $ and $H_n \in \diff^\infty(M, \lambda)$.

We choose a sequence $\alpha_n := p_n' / q_n'$ such that $|\alpha_n
- \alpha | \rightarrow 0$ monotonically. This choice defines a
sequence of intermediate scales by $q_n = q_{n-1}^d q_n'$ satisfying
$q_n' < q_n < q_{n+1}'$ which are geometrically natural for all the
previous transformations. Fixing $q_n$ determines $H_{n+1}$ via the
iterative formula
\begin{equation}
  \label{eq:HnDefn}
  H_{n+1} = H_n h_{n,q_n}.
\end{equation}
Defining the family of maps $h_{n,q}$ and investigating their
properties will form the bulk of this paper. 

\subsection{Reduction}

Though Theorem \ref{thm:GeneralCase} appears considerably more general
than Theorem \ref{thm:TorusCase} they follow from nearly identical
arguments. We are able to reduce the case of a general $M$ admitting a
smooth $C^\infty$ action of $\T$ to the case of $M = I^{d-1} \times
\T$, where $I=[0,1]$ is the standard unit interval, with $S_\theta:
I^{d-1} \times \T \rightarrow I^{d-1} \times \T$ given by $$S_\theta
(x_1, \dots , x_d) = (x_1, \dots ,x_{d-1}, x_d + \theta \bmod 1).$$

Let $\sigma$ denote the effective $\T$ action on $M$. For $q \geq 1$
we denote by $F_q$ the set of fixed points of the map
$\sigma(1/q,\cdot)$ and let $B:= \partial M \cup \bigcup_{q\geq1} F_q$
be the set of exceptional points.

We quote the following proposition of \cite{FK} that is similar to
other statements in \cite{AK,Windsor}

\begin{prop}\cite[proposition 5.2]{FK} \label{prop:reduction} Let $M$ be
  an $d$-dimensional compact connected $C^\infty$ manifold with an
  effective circle action $\sigma$ preserving a smooth volume $\mu$.
  Then here exists a continuous surjective map $\Gamma: I^{d-1}\times
  \T \to M$ with the following properties
\begin{enumerate}
\item The restriction of $\Gamma$ to $ (0,1)^{d-1} \times \T$ is a
  $C^{\infty}$ diffeomorphic embedding;
\item $\mu(\Gamma(\partial (I^{d-1}\times \T ))=0$;
\item $\Gamma(\partial(I^{d-1}\times \T)) \supset B$;
\item$\Gamma_*(\lambda)=\mu$;
\item $\sigma \circ \Gamma = \Gamma \circ S$.
\end{enumerate}
\end{prop}

An application of Proposition \ref{prop:reduction} at each step allows
us to conclude Theorem \ref{thm:GeneralCase} from the special case $M
= I^{d-1}\times \T$. Thus the construction need only be carried out
for two specific manifolds; $M = \T^d$ or $M = I^{d-1}\times \T$. For
both we take the action $S_\theta: M \rightarrow M$ given
by $$S_\theta (x_1, \dots , x_d) = (x_1, \dots ,x_{d-1}, x_d + \theta
\bmod 1)$$ that preserves the smooth unit volume $\lambda$ induced by
the usual Lebesgue measure on $\R^d$.

\subsection{Partitions and Measure-Theoretic Isomorphism}

The most difficult property to define on a finite scale is that of
measure-theoretic isomorphism to a circle rotation. We use the
abstract theory of Lebesgue spaces. Given an isomorphism of measures
space $(M_1, \mathfrak{B}_1,\mu_1 )$ and $(M_2, \mathfrak{B}_2,\mu_2
)$ there is a natural isomorphism of the associated
measure-algebras. If both the measure-spaces are Lebesgue spaces then
the converse is true; every isomorphism of the measure-algebras arises
from a point isomorphism of the measure spaces.  This is the crucial
observation that leads to the follwing abstract lemma, which appears
as \cite[Lemma 4.1] {AK}.

Given a partition $\xi$ of a space $M$ we write $\xi(x)$ for the atom of
the partition which contains $x$. We say that a sequence of partitions
$\xi_n$ generates if there is a set $F$ of full measure such that for
every $x \in F$ we have
\begin{equation*}
  \{ x \} = F \cap \bigcap_{n=1}^\infty \xi_n (x).
\end{equation*}

\begin{lem}\label{lem:Isomorphism}
  Let $M_1$ and $M_2$ be Lebesgue spaces. Let
  $(\xi_n^{(i)})_{n=1}^\infty$ be a monotone sequence of finite
  measurable partitions of $M_i$ that generates. Let
  $(T_n^{(i)})_{n=1}^\infty$ be a sequence of automorphisms of $M_i$ such
  that
  \begin{enumerate}
  \item $(T_n^{(i)})_{n=1}^\infty$ converges in the weak topology to
    an automorphism $T^{(i)}$ of $M_i$.
  \item $T_n^{(i)} \xi_n^{(i)} =  \xi_n^{(i)}$.
  \end{enumerate}
  Suppose that for each $n$ there exists a measure-theoretic
  isomorphism $K_n : M_1 / \xi_n^{(1)} \rightarrow M_2 / \xi_n^{(2)}$
  of the probability vectors such that:
   \begin{enumerate}
   \item $K_n^{-1} \circ \Bigl. T_n^{(2)} \Bigr|_{\xi_n^{(2)}} \circ
     K_n = \Bigl. T_n^{(1)} \Bigr|_{\xi_n^{(1)}}$.
   \item for all $\Delta \in \xi_{n-1}^{(1)}$ 
     $$K_n \Delta = K_{n-1} \Delta.$$
   \end{enumerate}
   Then the automorphisms $T^{(1)}$ and $T^{(2)}$ are
   measure-theoretically isomorphic.
\end{lem}


Consider the partition of $\T$ given by 
\begin{equation}
  \label{eq:partition1}
  \tilde \eta_{q} := \{ \tilde \Delta_{i,q} : 0 \leq i < q^d \}
\end{equation}
where $\tilde \Delta_{i,q} := [i q^{-d},(i+1)q^{-d}).$ This partition
is preserved under the action $R_{p/q}$. For any increasing sequence
of $q_n$ the sequence of partitions $\tilde \eta_{q_n}$ generates. Let
$M_2 = \T$, $\xi^{(2)}_n = \tilde \eta_{q_n}$ and $T_n^{(2)} =
R_{\alpha_n}$. Since $q_{n}$ divides $q_{n+1}$ we have $ \tilde\eta_{q_n} <
\tilde\eta_{q_{n+1}}$.

Let $\pi_d : M \rightarrow \T$ denote the projection onto the last
component of $M$. We obtain a partition of $M$ by
\begin{equation}
  \eta_{q} = \pi_d^{-1} \tilde \eta_{q} = \{ \Delta_{i,q} : 0 \leq i <
  q^d \}
\end{equation} where
\begin{equation*}
  \Delta_{i,q} := \{x : x_d \in [i q^{-d},(i+1)q^{-d}) \}.
\end{equation*}
Since $ \pi_d \circ S_\alpha = R_\alpha \circ \pi_d$ the partition
$\eta_q$ is preserved under the action of $S_{p/q}$ and, moreover, the
action of $S_{p/q}$ on $\eta_q$ is conjugated with that of $R_{p/q}$
on $\tilde \eta_{q}$.  Unfortunately the sequence of partitions
$\eta_{q_n}$ does not generate.

\begin{figure}[h]
  \psfrag{A}{$\pi_d$}
  \psfrag{1}[B]{$\Delta_{1,3}$}
  \psfrag{11}[Bc]{$\tilde \Delta_{1,3}$}
  \psfrag{2}[B]{$\Delta_{9,3}$}
  \psfrag{12}[Bc]{$\tilde \Delta_{9,3}$}
  \centering
  \includegraphics[scale=.8]{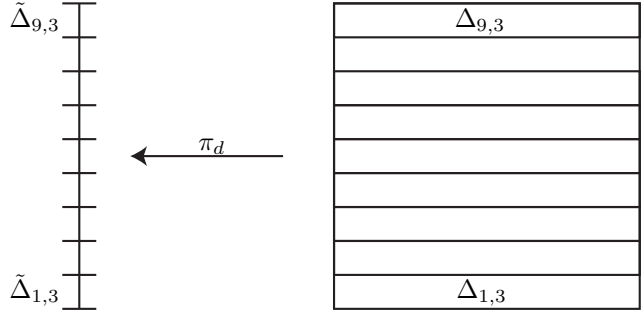}
  \caption{The partition $\eta_3$ of either $I\times\T$ or $\T^2$ and
    the partition $\tilde \eta_3$ of $\T$. }
\end{figure}
 
Let $M_1 = M$ and define the sequence of partitions 
\begin{equation}
  \label{eq:xinDefn}
  \xi_n^{(1)} := H_{n+1} \eta_{q_n} = H_n h_{n,q_n} \eta_{q_n}.  
\end{equation}
Unlike the sequence $\eta_{q_n}$, the sequence $\xi_n^{(1)}$ can be made to
generate. We construct $h_{n,q}$ as a diffeomorphism of $\pi_d^{-1}
[0, q^{-1}]$ and extend it to all of $M$ by requiring that it commute
with $S_{q^{-1}}$. Then 
\begin{enumerate}
\item Since $q_{n-1}^d$ divides $q_n$ we have for $0 \leq i < q_{n-1}^d$ 
  $$h_{n,q_n} \Delta_{i, q_{n-1}} = \Delta_{i,q_{n-1}}.$$

\item Since $q_n'$ divides $q_n$ we have 
  $$h_{n,q_n} \circ S_{\alpha_n} = S_{\alpha_n} \circ h_{n,q_n}.$$
\end{enumerate}
As $\eta_{q_{n-1}} < \eta_{q_n}$ we have  $H_{n+1}
\eta_{q_{n-1}} < H_{n+1} \eta_{q_n}$. By the first of our two
properties we have that $H_{n+1} \eta_{q_{n-1}} = H_{n}
\eta_{q_{n-1}}$ and hence $\xi_{n-1}^{(1)} < \xi_n^{(1)}$. Thus
$\{\xi_n^{(1)}\}$ is a monotone sequence of partitions as required
by Lemma \ref{lem:Isomorphism}. The second property ensures that $T_n
\xi_n^{(1)} = \xi_n^{(1)}$. Define the map 
\begin{equation*}
  K_n = \pi_d \circ H_{n+1}^{-1}.
\end{equation*}
Using the two properties we have that
\begin{align*}
  K_n \circ T_n^{(1)} &= T_n^{(2)} \circ K_n\\
  K_n (H_n \Delta_{i,q_{n-1}}) &= K_{n-1} (H_n \Delta_{i,q_{n-1}})
\end{align*}
as required by Lemma \ref{lem:Isomorphism}.

This completes the proof of the main theorem except for the proof that
the sequence $T_n$ converges in $\diff^\infty(M, \lambda)$ and the
proof that $\xi_{n}^{(1)}$ generates. 

\subsection{Construction of the Conjugating Maps.}

We will carry out the constructions for $M = \T^d$ and $M= I^{d-1}
\times \T$ simultaneously. The proof of unique ergodicity in the case
$M=\T^d$ will appear in a later section. 

\begin{lem}\label{lem:hndefn}
  Let  $n > 2 d$ and $q \in N$. There exists a map $h_{n,q} \in
  \diff^\infty( M , \lambda ) $ and a set $E_{n,q} \subset M$ such
  that:
  \begin{enumerate}
  \item $h_{n,q} S_{q^{-1}} = S_{q^{-1}} h_{n,q}$ and $h_{n,q}\bigl(
    \pi_d^{-1}[0,q^{-1}]\bigr) = \pi_d^{-1}[0,q^{-1}]$.

  \item $\lambda ( E_{n,q} ) >  1- 4 \frac{d -1}{n^2} $.
    
  \item for each $0\leq i < q^d$, 
    $$\diam h_{n,q} ( \Delta_{i,q} \cap E_{n,q} ) < \sqrt{d}
    q^{-1}.$$
    
  \end{enumerate}
\end{lem}

\subsubsection{Heuristic Construction.}

In order to motivate the construction of the family of conjugacy maps
we first construct a family of measure-preserving discontinuous maps
$\tilde h_{q}$ such that $\tilde h_q$ commutes with $S_{q^{-1}}$ and
carries each $\Delta_{i,q}$ into a $d$-dimensional cube with
side-length $q^{-1}$. 

\begin{figure}[h]
  \centering
  \psfrag{1}[Bc]{$\Delta_{1,3}$}
  \psfrag{2}[Bc]{$\Delta_{2,3}$}
  \psfrag{3}[Bc]{$\Delta_{3,3}$}
  \psfrag{4}[Bc]{$\Delta_{4,3}$}
  \psfrag{5}[Bc]{$\Delta_{5,3}$}
  \psfrag{6}[Bc]{$\Delta_{6,3}$}
  \psfrag{7}[Bc]{$\Delta_{7,3}$}
  \psfrag{8}[Bc]{$\Delta_{8,3}$}
  \psfrag{9}[Bc]{$\Delta_{9,3}$}
  \psfrag{11}[Bc]{$\tilde \phi_3 \Delta_{1,3}$}
  \psfrag{12}[Bc]{$\tilde \phi_3 \Delta_{2,3}$}
  \psfrag{13}[Bc]{$\tilde \phi_3 \Delta_{3,3}$}
  \psfrag{14}[Bc]{$\tilde \phi_3 \Delta_{4,3}$}
  \psfrag{15}[Bc]{$\tilde \phi_3 \Delta_{5,3}$}
  \psfrag{16}[Bc]{$\tilde \phi_3 \Delta_{6,3}$}
  \psfrag{17}[Bc]{$\tilde \phi_3 \Delta_{7,3}$}
  \psfrag{18}[Bc]{$\tilde \phi_3 \Delta_{8,3}$}
  \psfrag{19}[Bc]{$\tilde \phi_3 \Delta_{9,3}$}
  \psfrag{A}[Bc]{$\tilde \phi_3$}
  \includegraphics[scale=.8]{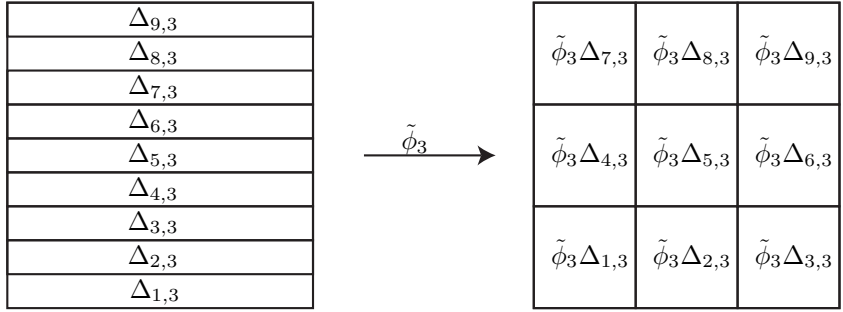}
  \caption{Action of $\tilde \phi_3 = \tilde h_3$ on the partition
    $\eta_3$.}
\end{figure}
   
Let $\tilde \phi_q $ be defined on $[0,1] \times[0, q^{-1}]$ by
letting it act on the interior by 
$$
\tilde \phi_q (x,y) := (q y, q^{-1} (1-x)  )
$$
and extend it to all of $[0,1] \times [0,1]$ by requiring $\tilde
\phi_q (x,y+q^{-1}) = \tilde \phi_q (x,y) + (0, q^{-1})$. Define
$\tilde \phi_{q}^{(i)}$ by
\begin{equation}
  \label{eq:phiqDefn}
  [ \tilde \phi_{q}^{(i)}]_j (x_1, \dots , x_d) = 
  \begin{cases}
    [\tilde\phi_q]_1 (x_i,x_{i+1}) & j = i\\
    [\tilde\phi_q]_2 (x_i,x_{i+1}) & j = i+1\\
    x_j & \text{otherwise}
  \end{cases}
\end{equation}
The map $\tilde h_{q}$ is defined by
\begin{equation*}
\tilde h_{q} := \tilde \phi_{q}^{(1)} \cdots \tilde
\phi_{q}^{(d-1)}.  
\end{equation*}
Each $\Delta_{i,q}$ is mapped, by $\tilde h_{q}$, into a cube of
side-length $q^{-1}$. The map $\tilde h_{q}$ commutes with
$S_{q^{-1}}$ since $\tilde\phi_{q}^{(d-1)}$ commutes with $S_{q^{-1}}$
by construction and the other $\tilde\phi_{q}^{(i)}$ don't affect
$x_d$.

\subsubsection{Proof of Lemma \ref{lem:hndefn}}

Our family of conjugating maps $h_{n,q}$ is constructed using the same
process as $\tilde h_{q}$ above.  Clearly control of some of the space
must be relinquished in order to be able to produce a $C^\infty$
volume preserving map. One additional complication arises ensuring
that we retain sufficient control over every orbit. Let $\varphi_n$
denote a $C^\infty$ map of the unit square satisfying
\begin{enumerate}
\item $\varphi_n$ is $C^\infty$ flat on the boundary.
  
\item $\varphi_n$ acts as a pure rotation by $\frac{\pi}{2}$ on $
  \bigl[ \frac{1}{n^2} , 1-\frac{1}{n^2} \bigr] \times \bigl[ \frac{1}{n^2}
  , 1-\frac{1}{n^2} \bigr].$
  
\item $\varphi_n$ preserves Lebesgue measure.
\end{enumerate}
Let $C_q(x,y) := (x, q^{-1} y)$ and define $\phi_{n,q}$ on $[0,1]
\times [0 , q^{-1}]$ by
\begin{equation}
  \label{eq:phinqDefn}
  \phi_{n,q} := C_q \varphi_n C_q^{-1}.
\end{equation}
Extend $\phi_{n,q}$ to the entire unit square by requiring that 
$$\phi_{n,q}(x,y + q^{-1} ) = \phi_{n,q}(x,y) + (0, q^{-1}).$$ 
This agrees with $\tilde \phi_{n,q}$ on a set of volume $(1-2/n^2)^2$
which we estimate from below by $1-4/n^2$. Analogously to our earlier
definition of $\tilde \phi_q^{(i)}$ we define $\phi_{n,q}^{(i)}$.
\begin{equation*}
  [ \phi_{q}^{(i)}]_j (x_1, \dots , x_d) = 
  \begin{cases}
    [\phi_q]_1 (x_i,x_{i+1}) & j = i\\
    [\phi_q]_2 (x_i,x_{i+1}) & j = i+1\\
    x_j & \text{otherwise}
  \end{cases}
\end{equation*}

\subsubsection{$M = I^{d-1} \times \T$ Case.} 
We define the conjugating map $h_{n,q}:I^{d-1}\times \T \rightarrow
I^{d-1}\times \T$ by
\begin{equation*}
  h_{n,q}:= \phi_{n,q}^{(1)} \cdots \phi_{n,q}^{(d-1)}.
\end{equation*}
 This map agrees with $\tilde h_q$ on a set $E_{n,q}$ given by 
\begin{equation}
  \label{eq:EnqDefinition}
  E_{n,q}^c = \bigcup_{i=1}^{d-1} \pi_i^{-1}\Bigr(
  \bigl[0,\frac{1}{n^2}\bigr) \cup \bigl(1-\frac{1}{n^2}, 1 \bigr]\Bigr)
  \cup \bigcup_{j=1}^{d-1} \bigcup_{k=1}^{q^j} 
  \pi_d^{-1} ( \frac{k}{q^j} - \frac{1}{n^2 q^j} , \frac{k}{q^j} +
  \frac{1}{n^2 q^j}). 
\end{equation}
Treating the sets on the right as disjoint we can estimate
\begin{equation}\label{eq:En}
  \lambda(E_{n,q}) > 1- 4 \frac{d-1}{n^2}.
\end{equation}

\subsubsection{$M = \T^d$ Case.}

In order to produce a unique ergodic diffeomorphism $T$ it is
necessary to control \emph{all} orbits. The set $E_{n,q}$ constructed
above for the case of $M = I^{d-1} \times \T$ excludes entire
orbits. In order to rectify this requires one more map. Let $\psi_q:
\T^d \rightarrow \T^d$ denote the translation
\begin{equation}
  \label{eq:psiqDefn}
  \psi_{q} (x_1, \dots, x_{d-1},x_d) := (x_1, \dots, x_{d-1},x_d) +
  x_d (q, \dots, q, 0 ) \mod 1.
\end{equation}
Obviously $\psi_q$ commutes with $S_{q^{-1}}$ and preserves the
Lebesgue measure. Furthermore, since $\psi_q$ does not affect the last
coordinate, it preserves each $\Delta_{i,q}$. For the uniquely
ergodic case we define 
\begin{equation}
  \label{eq:hn}
  h_{n,q} := \phi_{n,q}^{(1)} \cdots \phi_{n,q}^{(d-1)} \psi_q
\end{equation}
Exactly as for the ergodic case $h_{n,q}$ agrees with $\tilde h_q$ on a
set $E_{n,q}$ with
$$\lambda(E_{n,q}) > 1- 4 \frac{d-1}{n^2}.$$
The map $\psi_q$ ensures that $E_{n,q}$ contains most of \emph{every}
orbit.

\begin{figure}
  \centering
  \includegraphics[scale=.8]{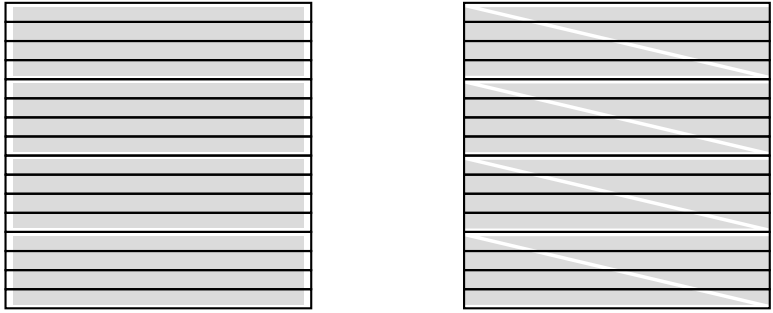}
  \label{fig:Enq}
  \caption{The set $E_{n_q}$ for the case $M = I \times \T$ (left) and
    for the case $M=\T^2$ (right).}
\end{figure}

\newpage
\subsection{Analytic Properties}

\subsubsection{Notation}\label{sec-analytic}

All of our diffeomorphisms $h : I^{d-1} \times \T \rightarrow I^{d-1}
\times \T$ are identity in a neighborhood of the boundary and hence
can be identified with a diffeomorphism $h : \T^d \rightarrow
\T^d$. Defining a topology on $\diff^k(\T^d,\T^d)$ defines a topology
on the closure of the space of diffeomorphisms $h : I^{d-1} \times \T
\rightarrow I^{d-1} \times \T$ that are identity in a neighborhood of
the boundary.

Let $f,g \in C^0(\T^d,\T^d)$. We define 
\begin{equation*}
  \hat{d}_0(f,g) = \max_{x \in M} d\bigl(f(x),g(x) \bigr).
\end{equation*}

Let $f \in C^k(\R^d,\R)$. Given $a \in \N^d$ we denote $|a| := a_1 +
\dots + a_d$ and
\begin{equation*}
  D_a f  := \frac{\partial^{|a|} f}{\partial x_1^{a_1} \dots \partial
  x_d^{a_d} }.
\end{equation*}
Using this we can define 
\begin{equation*}
  \tnorm{f}_k = \max_{1 \leq |a| \leq k} \max_{x \in M} |D_a f (x)|.
\end{equation*}
For $f \in C^k(\R^d,\R^d)$ we define
\begin{equation*}
  \tnorm{f}_k = \max_{1 \leq i \leq d}  \max_{1 \leq |a| \leq k}
  \max_{x \in M} |D_a f_i (x)|.
\end{equation*}
For $h : \T^d \rightarrow \T^d$ we can define a natural lift $\hat{h}:
\R^d \rightarrow \R^d$. Now given $f,g \in C^k (\T^d,\T^d)$ we define
\begin{equation*}
  \hat{d}_k( f, g) = \max\{ d_0(f,g), \tnorm{\hat{f} - \hat{g}}_k \} 
\end{equation*}

Finally, for $f,g \in \diff^k (\T^d,\T^d)$ we define
\begin{equation*}
  d_k(f,g) = \max \{ \hat{d}_k(f,g) , \hat{d}_k( f^{-1},g^{-1}) \} 
\end{equation*}

The metric defined in this way is equivalent to the usual one defined
via the operator norms but is easier to work with for explicit
estimates. For further details consult \cite{Saprykina}.

\subsubsection{Estimates}
\begin{lem}
  We have the following estimate:
  \begin{equation}
    \label{eq:hnqEstimate}
    |\!|\!| h_{n,q} |\!|\!|_k < C_1 q^{d k}
  \end{equation}
  where $C_1$ depends on $d$, $k$, and $n$ but is independent of $q$. 
\end{lem}

\begin{proof}
  By direct computation we obtain
  \begin{equation}
    \label{eq:estimate1}
    |\!|\!| \phi_{n,q}^{(i)} |\!|\!|_k < q^k |\!|\!| \varphi_n |\!|\!|_k
  \end{equation}
  and
  \begin{equation}
    \label{eq:estimate2}
     |\!|\!| \psi_q |\!|\!|_k < q.
  \end{equation}
  We claim that partial derivatives with $|a|=k$ consist of sums of
  products of at most $(d-1) k$ terms of the form
  \begin{equation}
    \label{eq:term3}
    \bigl( D_b [\phi_{n,q}^{(i)}]_j\bigr) ( \phi_{n,q}^{(i+1)} \dots
    \phi_{n,q}^{(d-1)} \psi_q )
  \end{equation}
  with $|b| \leq k $ and at most $k$ terms of the form
  \begin{equation}
    \label{eq:term4}
    D_c [\psi_q]_j
  \end{equation}
  with $|c|=1$. This is true for $|a|=1$ by computation and, by the
  product and chain rules, if it is true for $|a|=k$ then it is true
  for $|a|=k+1$. By induction it is therefore true for all $k$.
  
  Now suppose the estimate \eqref{eq:hnqEstimate} holds for $k$ we
  wish to show it holds for $k+1$. We use our structure theorem for
  $k$. Differentiating a term of the from \eqref{eq:term3} we get a
  sum of products of $d + 1- i$ terms. The first is of the form
  \eqref{eq:term3} but with the power of the derivative raised by $1$.
  The next $d-1-i$ terms are first partial derivatives of
  $\phi_{n,q}^{(i+1)}, \dots, \phi_{n,q}^{(d-1)}$. The final
  term is a first partial derivative of $\psi_q$. Applying the estimates
  \eqref{eq:estimate1} we see that the required power of $q$ has been
  increased by at most $d$. Differentiating \eqref{eq:term4} gives
  zero since $\psi_q$ is linear.
 
\end{proof}

By an application of the Fa\`a di Bruno's formula we obtain the
following corollary. 
\begin{cor}
  We have the following estimate
  \begin{equation}
    \label{eq:HnEstimate}
    |\!|\!| H_n h_{n,q} |\!|\!|_k < C_2 q^{k d} 
  \end{equation}
  where $C_2$ depends on $H_n$, $n$, and $k$ but is independent of $q$.  
\end{cor} 

\subsection{Completing the Construction}

Having now constructed the family of maps $h_{n,q}$ from which the
maps $H_n$ are assembled it remains only to explain how we choose the
sequence $q_n$. The choice of $q_n$ determines $\alpha_n$ as the best
approximation to $\alpha$ with denominator $q_n$.The choices of
$q_1,...,q_{n-1}$ completely determines $H_n$. We show how given $H_n$
we choose $q_n$ so that $T_n$ has the desired properties.

In the original Anosov and Katok method of construction the choice of
$\alpha_n$ in the definition of $T_n$ \eqref{eq:TnDefn} determined the
distance between the already determined $T_{n-1}$ and $T_n$ in
$\diff^n$. The observation there was that if $\alpha_n$ could be
chosen arbitrarily close to $\alpha_{n-1}$ then the transformation
$T_n$ could be made arbitrarily close to $T_{n-1}$. The advantages of
this approach are that no estimates on the maps $H_n$ are
required. Unfortunately this approach is inconsistent with ensuring
that the sequence $\alpha_n$ converges to an a priori given number
$\alpha$. In the approach we take the choice of $q_n$ (and hence of
$\alpha_n$) determines the distance between $T_n$ and, the as-yet
undetermined transformation, $T_{n+1}$. Since the choice of $q_n$
fixes the conjugacy map $H_{n+1}$ the only undetermined quantity in
$T_{n+1}$ is the choice of $\alpha_{n+1}$. Supposing only that the
choice of $\alpha_{n+1}$ will be a better approximation to $\alpha$
than $\alpha_n$ we are able to estimate the distance between $T_n$ and
$T_{n+1}$ knowing only the choice of $\alpha_n$.

\begin{lem}\label{lem:MainLem}
  Let $k \in \N$. For all $h \in \diff^k(M)$ and all $\alpha, \beta
  \in \R$ we obtain
  \begin{equation*}
    d_k(h \circ S_{\alpha}\circ h^{-1}, h \circ S_{\beta} \circ h^{-1})
    \leq C_3  \tnorm{h}_{k+1}^{k+1} | \alpha - \beta|
  \end{equation*}
  where $C_3$ depends only on $k$.   
\end{lem}

\begin{proof}
  For $k=0$ we have the estimate
  $$ d_0(h \circ S_{\alpha}\circ h^{-1}, h \circ S_{\beta} \circ h^{-1})
  \leq \tnorm{h}_1 | \alpha - \beta|$$ by the mean value theorem. We
  claim that for $a \in \N^d$ with $|a|=k$ the partial derivative 
  $$
  D_a [  h_i\circ S_{\alpha} \circ h^{-1} -  h_i\circ
  S_{\beta} \circ h^{-1} ]
  $$ 
  will consist of a sum of terms with each term being the product of a
  single partial derivative
  \begin{equation}
    \label{eq:term1}
    \bigl( D_b h_i \bigr) (S_{\alpha} h^{-1} ) -\bigl( D_b  h_i\bigr)
    (S_{\beta} h^{-1})
  \end{equation}
  with $|b| \leq k$, and at most $k$ partial derivatives of the form 
  \begin{equation}
    \label{eq:term2}
    D_b h^{-1}_j 
  \end{equation}
  with $|b| \leq k$. For $k=1$ we have 
  \begin{multline*}
    \frac{\partial}{\partial x_j} [  h_i\circ S_{\alpha} \circ h^{-1} -  h_i\circ
    S_{\beta} \circ h^{-1}] \\
    = \sum_{l=1}^d \bigl(\frac{\partial h_i
    }{\partial x_l} \circ S_{\alpha}\circ h^{-1} - \frac{\partial h_i
    }{\partial x_l} \circ S_{\beta} \circ h^{-1} \bigr) \frac{\partial
      h^{-1}_l}{\partial x_j}.
  \end{multline*}
  We proceed by induction. By the product rule we need only consider
  the effect of differentiating \eqref{eq:term1} and \eqref{eq:term2}.
  Differentiating \eqref{eq:term1} with respect to $x_j$ we obtain
  \begin{equation*}
    \sum_{l=1}^d \bigl(  \frac{\partial  D_b h_i}{\partial x_l} \circ
    S_{\alpha} \circ h^{-1} 
    - \frac{\partial  D_b h_i}{\partial x_l} \circ S_{\beta} \circ h^{-1}
    \bigr)\frac{\partial h^{-1}_l}{\partial x_j}.
  \end{equation*}
  which increases the number of terms of the form \eqref{eq:term2} by
  1.  Differentiating \eqref{eq:term2} we get another term of the form
  \eqref{eq:term2} but with $|b| \leq k+1$. 
  
  We estimate
  \begin{gather*}
    \| D_a h_i \circ S_{\alpha} \circ h^{-1} - D_a h_i \circ S_{\beta}
    \circ h^{-1} \|_0 \leq  \tnorm{h}_{|a|+1} | \alpha - \beta |\\
    \| D_a h^{-1}_l \|_0 \leq \tnorm{h}_{|a|}
  \end{gather*}
  These estimates together with claimed structure of the
  partial derivatives, and the fact that the inverse maps have the
  same structure, completes the proof. The constant $C_3$ is the
  number of terms in the sum which depends only on $k$ and not on the
  map $h$.
\end{proof}

Define $ F_n := H_{n+1} ( E_{n,q_n}) $ and let $F := \liminf
F_n$. Clearly, from Lemma \ref{lem:hndefn}, we have that 
\begin{equation*}
  \lambda(F) \geq \lim_{n \rightarrow \infty} (1 - 4 (d-1)
  \sum_{m=n}^\infty \frac{1}{m^2}) = 1.
\end{equation*}
We will show that any point in $F$ has a unique coding relative to the
sequence of partitions $\xi_n$.

\begin{prop}\label{prop:Completion}
  Let $\epsilon_n$ be a summable sequence of positive numbers. There
  is a choice of $\{q_n'\}$ such the transformations $T_n$ defined by
  \eqref{eq:TnDefn} satisfy
  \begin{enumerate}
  \item $d_n ( T_n , T_{n+1} ) < \epsilon_n$.

  \item for $A \in \xi_n $  
    $$\diam (A \cap F_n) < \epsilon_n $$  
  \end{enumerate}
\end{prop}

\begin{proof}
  By the definition of a Liouville number for any polynomial $P(q')$ we can
  find $q'_n> q_{n-1}$ such that $\alpha_n := p'_n / q'_n$ is a better
  approximation to $\alpha$ than $\alpha_{n-1}$ and such 
  \begin{equation*}
    P(q'_n) \bigl| \frac{p_n'}{q_n'} - \alpha \bigr| < \epsilon_n
  \end{equation*}
  We will define $q = q_{n-1}^d q'$ to ensure that $h_{n,q}$
  satisifes . Since $q < (q')^{d+1}$ we that for any polynomial $P(q)$
  we can find find $q'_n$ such that $\alpha_n := p'_n / q'_n$ is a better
  approximation to $\alpha$ than $\alpha_{n-1}$ and such 
  \begin{equation*}
    P(q_n) \bigl| \frac{p_n'}{q_n'} - \alpha \bigr| < \epsilon_n
  \end{equation*}
  Now combining \eqref{eq:HnEstimate} and Lemma \ref{lem:MainLem} we have 
  \begin{align*}
    d_{n} (T_n,T_{n+1}) &< P(q_n) | \alpha_n - \alpha_{n+1} |\\
    &< 2 P(q_n) | \alpha_n - \alpha |.
  \end{align*}
  Similarly for $H_{n+1} \Delta_{i,q_n} \in \xi_n$ we have 
  \begin{align*}
    \diam ( H_{n+1}\Delta_{i,q_n} \cap F_n) &=  \diam ( H_{n}
    h_{n,q_n} (\Delta_{i,q_n} \cap E_{n,q_n}))\\
    &\leq \| H_n \|_1 \diam h_{n,q_n} (\Delta_{i,q_n} \cap E_{n,q_n})\\
    &\leq \| H_n \|_1 \sqrt{d} q_n^{-1} 
  \end{align*}
  using Lemma \ref{lem:hndefn}.  Thus we see that we can choose
  $\alpha_n$ such that the required two properties hold.
\end{proof}

Since $\epsilon_n$ is summable we have that $\{T_n\}$ is a Cauchy
sequence in $\diff^\infty(M, \lambda)$ and hence converges to some
$T\in \diff^\infty(M, \lambda)$. For any $x \in F$ we have $x \in F_n$
for all but finitely many $n$. Thus, by Proposition
\ref{prop:Completion}, we have for all $x \in F$
$$\bigcap_{n=1}^\infty \xi_n(x) \cap F = \{ x\}.$$ 
This shows that $ \{ \xi_n\} $ is a generating partition and hence
completes the proof of Theorem \ref{thm:GeneralCase}. 
\section{Unique Ergodicity}

When $M = \T^d$ we wish to prove unique ergodicity. We will use the
following abstract lemma, also used in \cite{Windsor}.

\begin{lem}\label{lem:UniqueErg}
  Let $q_n$ be an increasing  sequence of natural numbers and $T_n : X
  \rightarrow X$ a sequence of transformations which converge
  uniformly to a transformation $T$. Suppose that for each
  continuous function $\varphi$ from a dense set of continuous
  functions $\Phi$ there is a constant $c$ such that
  \begin{equation}
    \label{eq:AnDefn}
    \frac{1}{q_n} \sum_{i=0}^{q_n -1} \varphi (T_n^i x)
    \xrightarrow[n \rightarrow \infty]{} c
    \text{ uniformly} 
  \end{equation}
  and
  \begin{equation}
    \label{eq:dqnDefn}
    d^{(q_n)} ( T_n , T) := \max_x \max_{0 \leq i < q_n} d (T_n^i x, T^i
    x ) \rightarrow 0
  \end{equation}
  Then $T$ is uniquely ergodic
\end{lem}
\begin{proof}
  Condition \eqref{eq:dqnDefn} implies that
  $$\| \frac{1}{q_n} \sum_{i=0}^{q_n -1} \varphi (T_n x) - \frac{1}{q_n}
  \sum_{i=0}^{q_n -1} \varphi (T x) \|_0 \rightarrow 0 $$ and then
  condition \eqref{eq:AnDefn} becomes the standard result that if the
  Birkhoff sums converge uniformly then the map is uniquely
  ergodic \cite{HK}. 
\end{proof}

To establish condition \eqref{eq:AnDefn} it is insufficient to know
only that $E_{n,q}$ has large measure, we also need to know that most
of every $S_\theta$ orbit intersects $E_{n,q}$.

For each $x \in \T^d$ define $\sigma_x : \T \rightarrow \T^d$ by
$\sigma_x \theta = S_\theta x$.

\begin{lem}
  Let $q > d n^2$.  For each $x \in \T^d$ there is a set
  $J_{n,q}^{(x)}\subset \T^d$, measurable with respect to $\eta_{q}$,
  with measure
  \begin{equation}
    \label{eq:Jn}
    \lambda(J_{n,q}^{(x)}) > 1- \frac{4 d}{n^2}
  \end{equation}
  such that if $\Delta_{i,q} \subset J_{n,q}^{(x)}$ then
  \begin{gather}
    \label{eq:GoodElement1}
      \sigma_x^{-1}( \Delta_{i,q} \cap E_{n,q}^c ) = \varnothing,\\
      \label{eq:GoodElement2}
      \lambda(\Delta_{i,q} \cap E_{n,q} ) > \bigl( 1- \frac{ 2 (d-1)
      }{n^2} \bigr) \lambda ( \Delta_{i,q} ).
  \end{gather}
\end{lem}

\begin{figure}
  \centering
  \psfrag{1}[Bc]{$\sigma_x(\T)$} 
  \includegraphics[scale=.8]{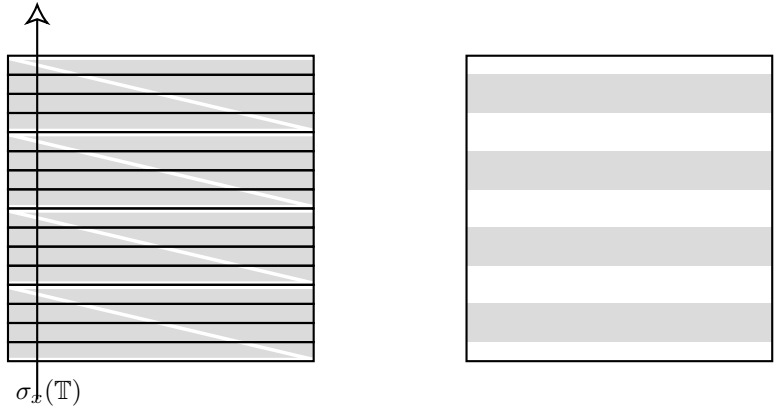}
  
  \caption{The orbit of $x \in \T^2$, indicated by the arrow on the
    left, combines with $E_{n,q}$, indicated by the shaded region on
    the left, to produce the set $J_{n,q}^{(x)}$, indicated by the
    shaded region on the right.}
  
\end{figure}

\begin{proof}
  It is immediate that
  \begin{equation}
    \label{eq:Enqc}
    (E_{n,q}')^c = \bigcup_{i=1}^{d-1} \pi_i^{-1} \bigl( - \frac{1}{n^2}
    , \frac{1}{n^2} ) \cup \bigcup_{j=1}^{d-1} \bigcup_{k=1}^{q^j}
    \pi_d^{-1} ( \frac{k}{q^j} - \frac{1}{n^2 q^j} , \frac{k}{q^j} +
    \frac{1}{n^2 q^j}) 
  \end{equation}
  Let $x$ be arbitrary. We compute $\sigma_x^{-1} \psi_q (E_{n,q}')^c$
  using \eqref{eq:Enqc} and \eqref{eq:psiqDefn}.
  \begin{align*}
    \sigma_x^{-1} \psi_q^{-1} \pi_i^{-1} \bigl( - \frac{1}{n^2} ,
    \frac{1}{n^2} \bigr) = \bigcup_{l=1}^{q} \Bigl( \frac{l}{q} -
    \frac{1}{n^2 q} - x_d - \frac{x_i}{q} , \frac{l}{q} + \frac{1}{n^2 q}
    - x_d - \frac{x_i}{q}\Bigr) \\
    \sigma_x^{-1} \psi_q^{-1} \pi_d^{-1} \bigl( \frac{k}{q^j} -
    \frac{1}{n^2 q^j} , \frac{k}{q^j} + \frac{1}{n^2 q^j} \bigr) = \bigl(
    \frac{k}{q^j} - \frac{1}{n^2 q^j} - x_d , \frac{k}{q^j} + \frac{1}{n^2
      q^j} - x_d \bigr)
  \end{align*}

  This excluded set of $\tau$ consists of at most $(d-1) q + q^{d-1}$
  intervals. Expanding these intervals to make them measurable with
  respect to $\sigma_x^{-1} \eta_{q}$ excludes an additional set of
  measure at most
  $$ 
  \frac{2}{q^d} \bigl( (d-1) q + q^{d-1} \bigr) < \frac{4}{n^2}.
  $$
  Let $E$ denote the measurable hull of $\sigma_x^{-1} E_{n,q}^c$
  in $\sigma_x^{-1} \eta_{q}$. We have $\lambda(E) = 4 d / n^2$.
  Define the set $J_{n,q}^{(x)}$ to be the $\eta_{q}$ measurable set
  satisfying
  $$
  \sigma_x^{-1} J_{n,q}^{(x)} = E^{c}.
  $$ 
\end{proof}
Note that the proportion in \eqref{eq:Jn} is lower than the proportion
in \eqref{eq:En}. We have had to give up control over parts of each
orbit in order to gain control over all orbits. The set
$J_{n,q}^{(x)}$ consists of those atoms of $\eta_q$ where we have
control over the behaviour of all of $S_\theta x$ under $h_{n,q}$.

Using the geometric information contained in these lemmas we can prove
a distribution result.  

\begin{prop}\label{prop:Distribution}
  Let $\epsilon >0$, $q \in \N$, and $\varphi$ be a $(\sqrt{d}
  q^{-d},\epsilon)$-uniformly continuous function, i.e 
  \begin{equation*}
    \varphi( B_{\sqrt{d} q^{-d}}(x)) \subset B_{\epsilon}(\varphi(x)).
  \end{equation*}
  For all $q' \in \N$ and for all $x \in \T^d$,
  \begin{equation}
    \label{eq:LocalEstimate}
    \biggl| \frac{1}{q'} \sum_{i=0}^{q'-1} \varphi
    ( h_{n,q} S_{1/q'}^{i} x ) - \int \varphi d\lambda \biggr| <
    \frac{14 d}{n^2} \| \varphi \|_{0} + \frac{2 q^d}{q'}  \| \varphi
    \|_{0} + 2 \epsilon. 
  \end{equation}
\end{prop}
\begin{proof}
  For $x,y \in \Delta_{i,q} \cap E_{n,q}$ we have 
  \begin{equation*}
    d( h_{n,q} x, h_{n,q} y) \leq \diam h_{n,q}(\Delta_{i,q} \cap E_{n,q}) \leq
    \sqrt{d} q^{-d}.
  \end{equation*}
  By the hypothesis on $\varphi$ we have $|\varphi(h_{n,q} x) -
  \varphi(h_{n,q}y)|<2 \epsilon$. Averaging over all $y \in
  \Delta_{i,q} \cap E_{n,q}$ we obtain for any  $x \in
  \Delta_{i,q} \cap E_{n,q}$,
  \begin{equation}
    \label{eq:LocalAverage}
    \biggl|\varphi(h_{n,q} x) - \frac{1}{\lambda(\Delta_{i,q} \cap E_{n,q}
      )} \int_{h_{n,q} (\Delta_{i,q} \cap E_{n,q})} \varphi\,
    d\lambda\biggr| < 2 \epsilon.
  \end{equation}

  Let $\mathcal{O}^{(x)}$ consist of $\lfloor \frac{q'}{q^d} \rfloor
  q^d$ points of the orbit of $x$ under $S_{1/q'}$ that are
  equidistributed among the atoms of the partition $\eta_q$. There are
  at most $q^d$ exceptional points outside of $\mathcal{O}^{(x)}$.

  By \eqref{eq:GoodElement1} for $\Delta_{i,q} \subset J_{n,q}^{(x)}$
  the number of points from $\mathcal{O}^{(x)}$ in $\Delta_{i,q} \cap
  E_{n,q}$ is $\lfloor\frac{q'}{q^d}\rfloor$.  Let $I := \{ 0\leq i <
  q' : S_{1/q'}^i x \in J_{n,q}^{(x)} \cap \mathcal{O}^{(x)} \}$ be
  the equidistributed points in good atoms. Using this count and
  \eqref{eq:LocalAverage} we obtain
  \begin{equation*}
    \begin{split}
      \biggl| \frac{1}{q'} \sum_{i \in I} \varphi ( h_{n,q} &S^i_{1/q'}
      x
      )\\
      &- \frac{1}{q'} \sum_{\Delta_{i,q} \subset J_{n,q}^{(x)}}
      \biggl\lfloor\frac{q'}{q^d}\biggr\rfloor
      \frac{1}{\lambda(\Delta_{i,q} \cap E_{n,q} )} \int_{h_{n,q}
        (\Delta_{i,q} \cap E_{n,q}) } \varphi\, d\lambda \biggr| < 2
      \epsilon.
    \end{split}
  \end{equation*}
  The remaining estimates just formalize the observation that since
  $J_{n,q}^{(x)}$ is nearly full measure and since $I$ is nearly all
  of the orbit the above estimate implies \eqref{eq:LocalEstimate}. 

  First we produce estimates that account for the fact that $q^d$ does
  not divide $q'$ and hence we do not have equidistribution of the
  entire orbit.
  \begin{gather*}
    \begin{split}
         \biggl| \frac{1}{q'} \sum_{\Delta_{i,q} \subset J_{n,q}^{(x)}}&
         \biggl\lfloor\frac{q'}{q^d}\biggr\rfloor
         \frac{1}{\lambda(\Delta_{i,q} \cap E_{n,q} )} \int_{h_{n,q}
           (\Delta_{i,q} \cap E_{n,q}) } \varphi \, d\lambda \\
         &- \sum_{\Delta_{i,q} \subset J_{n,q}^{(x)}}
         \frac{1}{q^d}\frac{1}{\lambda(\Delta_{i,q} \cap E_{n,q} )}
         \int_{h_{n,q} (\Delta_{i,q} \cap E_{n,q}) } \varphi \, d\lambda
         \biggr| < \frac{q^d}{q'} \| \varphi \|_0
       \end{split}\\
        \biggl| \frac{1}{q'} \sum_{i=0}^{q'-1} \varphi( h_{n,q} S^i_{1/q'}
    x ) -  \frac{1}{q'} \sum_{i \in \mathcal{O}^{(x)}} \varphi( h_{n,q} S^i_{1/q'}
    x ) \biggr| < \frac{q^d}{q'}\| \varphi \|_0
  \end{gather*}
  Second we produce estimates using \eqref{eq:Jn} and \eqref{eq:GoodElement1}
  \begin{gather*}    
      \biggl| \frac{1}{q'} \sum_{i \in \mathcal{O}^{(x)}} \varphi( h_{n,q} S^i_{1/q'}
      x ) -  \frac{1}{q'} \sum_{i \in I} \varphi( h_{n,q} S^i_{1/q'}
      x ) \biggr| < \frac{4 d}{n^2}\| \varphi \|_0,\\
      \biggl| \int_{h_{n,q}  J_{n,q}^{(x)} } \varphi \, d\lambda -  \int
    \varphi \, d\lambda \biggr| < \frac{ 4d }{n^2} \|\varphi\|_0 
  \end{gather*}
  Finally we produce estimates using \eqref{eq:GoodElement2}
  \begin{gather*}
    \biggl|\int_{h_{n,q} (J_{n,q}^{(x)} \cap E_{n,q}) } \varphi \,
    d\lambda - \int_{h_{n,q} J_{n,q}^{(x)} } \varphi \, d\lambda
    \biggr| < \frac{2 (d-1)}{n^2}\|\varphi\|_0,\\
      \biggl| \frac{1}{q^d \lambda (\Delta_{i,q} \cap E_{n,q})} 
    \int_{h_{n,q}  (J_{n,q}^{(x)} \cap E_{n,q}) } \varphi \, d\lambda -
    \int_{h_{n,q}  (J_{n,q}^{(x)} \cap E_{n,q}) } \biggr|
    <\frac{4(d-1)}{n^2} \|\varphi\|_0. 
  \end{gather*} 
Combining these estimates gives us exactly \eqref{eq:LocalEstimate} as
required. 
\end{proof}

Let $\Phi = \{ \varphi_n \}$ be a set of Lipshitz functions that is
dense in $C^0(M, \R)$. Let $L_n$ be a Lipshitz constant for
$\varphi_1\circ H_n, \dots, \varphi_n\circ H_n$. At step $n$ we can
choose $q_n'$ so that $L_n \sqrt{d}q_n^{-1} < n^{-2}$ and $q_n' > n^2
q_n$. Then applying Proposition \ref{prop:Distribution} we see that
for $\varphi \in \{ \varphi_1, \cdots, \varphi_n \}$ we have
\begin{equation*}
  \biggl| \frac{1}{q_n'} \sum_{i=0}^{q'_n-1} \varphi
  ( T_{n+1}^{i} x ) - \int \varphi d\lambda \biggr| <
  \frac{17 d}{n^2} \| \varphi \|_{0}.
\end{equation*}
This establishes \eqref{eq:AnDefn} from Lemma \ref{lem:UniqueErg}. To
establish \eqref{eq:dqnDefn} from \ref{lem:UniqueErg} observe that 
\begin{align*}
  d^{(q_n)}(T_n, T_{n+1} ) &\leq \tnorm{H_{n+1}}_1 q_n |\alpha_n -
  \alpha_{n+1}\\
  &\leq P(q_n) |\alpha_n - \alpha|
\end{align*}
and hence we can choose $q_n'$ so that this is less than $1/n$. In
actual fact this estimate is weaker than those that arise in the proof
of Proposition  \ref{prop:Completion} and so is automatic. 

This verifies the hypotheses of Lemma \ref{lem:UniqueErg} and hence we
conclude that $T$ is uniquely ergodic.

\end{document}